\documentclass[12pt]{article}
\usepackage{amsthm}
\usepackage[utf8]{inputenc}
\usepackage{amsmath,amsfonts}
\usepackage{graphicx,epstopdf}
\usepackage{color,hyperref}
\usepackage{subcaption}
\usepackage{xfrac}
\usepackage{framed}
\definecolor{shadecolor}{rgb}{.9, .9, .9}

\usepackage{pgfplots}
\pgfplotsset{compat=newest}

\newcommand{\N}{\mathbf N}

\newcommand{\F}{\mathbf F}

\newcommand{\Q}{\mathbf Q}
\newcommand{\Z}{\mathbf Z}

\newcommand{\C}{\mathbb C}

\graphicspath{ {./figures/experimental/} }

\newtheorem{theorem}{Theorem}

\newtheorem{lemma}[theorem]{Lemma}
\newtheorem{proposition}[theorem]{Proposition}

\newcommand{\sherh}[1]{\fboxsep=0pt\setlength{\fboxrule}{1pt}
\begin{center}
   \fbox{\colorbox{green}{
         \begin{minipage}[t]{13cm}
            #1
         \end{minipage}
      }
   }
\end{center}}
\newcommand{\sherhh}[1]{\fboxsep=0pt\setlength{\fboxrule}{1pt}
\begin{center}
   \fbox{\colorbox{yellow}{
         \begin{minipage}[t]{13cm}
            #1
         \end{minipage}
      }
   }
\end{center}}

\newcommand{\sherhhh}[1]{\fboxsep=0pt\setlength{\fboxrule}{1pt}
\begin{center}
   \fbox{\colorbox{red}{
         \begin{minipage}[t]{13cm}
            #1
         \end{minipage}
      }
   }
\end{center}}

\renewcommand{\sherh}[1]{}\renewcommand{\sherhh}[1]{}\renewcommand{\sherhhh}[1]{}

\title{Composition laws on the Fricke surface and Markov triples}
\author{A. Muhammed Uluda\u{g}\footnote{
{Department of Mathematics, Galatasaray University,}
{\c{C}{\i}ra\u{g}an Cad. No. 36, 34357 Be\c{s}ikta\c{s}}
{\.{I}stanbul, Turkey}},  Esra \"{U}nal Yılmaz \footnote{
{Department of Mathematics, Bolu Abant Izzet Baysal University,}
{ G\"{o}lk\"{o}y 14100,}
{Bolu, Turkey}}}

\begin{document}
\maketitle


\begin{abstract}
We determine some composition laws related to the Fricke surface and the ``double'' Fricke surface. This latter surface admits the squares of Markov triples as its solutions.
\end{abstract}

{\small \paragraph{Keywords}
Fricke surface, Markov equation, Markov triples, group law on quadrics, Frobenius conjecture.

\section{Introduction}\label{intro}
It is well known that the positive integral solutions of the {\it Markov equation}
\begin{align}\label{mark}
x^2+y^2+z^2=3xyz
\end{align}
can be obtained from the solution $(1,1,1)$ by using the 
{\it Viète transformations}
$$
L:(x,y,z)\rightarrow (x,3xy-z, y) \mbox{ and } R:(x,y,z)\rightarrow (y,3yz-x, z),
$$
and the permutations of these solution triples \cite{aigner}. The cubic surface $\mathbf F$ defined over some field by (\ref{mark}) is called the {\it Fricke surface}. An integral solution of (\ref{mark}) is called a {\it Markov triple}; except $(0,0,0)$, these are obtained from positive integral solutions by a pair of sign changes of coordinates. The famous conjecture of Frobenius states that the largest component of a positive Markov triple determines uniquely the remaining components.

Suppose $\mathbf H$ is a cubic hypersurface in the projective space $\mathbb P^n(\C)$ and $P, Q$ are two generic points on $\mathbf H$. Then the line through $P$ and $Q$ intersects $\mathbf H$ at a third point $R$. Thus we may introduce the partial composition law $P\circ Q:=R$ on $\mathbf H$. In case $\mathbf H$ is an elliptic curve and $O$ is a fixed point of $\mathbf H$, then the product 
$P\star Q:=O\circ(P\circ Q)$ defines an abelian group law on 
$\mathbf H$. However, in case $\dim \mathbf H\geq 3$, this composition is not always well-defined and does not yield a group law since it is not associative \cite{manin}. By linearity, this composition satisfies 
\begin{align}\label{ident} 
\alpha \circ \beta=\gamma \iff \alpha \circ \gamma=\beta \iff \beta \circ \gamma=\iff \alpha
\end{align}
Our aim in this paper is to determine this composition law on the Fricke surface $\mathbf F$ in $\C^3$ defined by (\ref{mark}) and study it. 

Sections of $\mathbf F$ with hyperplanes parallel to the coordinate axis gives quadric curves. Any smooth quadric $Q$ can be endowed with a genuine commutative group law, defined as follows \cite{poormans}:  Let $O\in Q$ be a (rational) point, and $A, B\in Q$. Let $L$ be the line through $O$ which is parallel to the line through $A$ and $B$. Then $A\oplus B$ is defined to be the second intersection point of $L$ with $Q$. We also determine this group law for these quadric hyperplane sections. The Frobenius conjecture on the Markov numbers can be formulated as a statement about the integral points on these quadric sections. 

Now consider the equation
\begin{align}\label{gen}
(x+y+z)^2=9xyz.
\end{align}
We call the surface defined by this equation the {\it double Fricke surface} and denote it by $\mathbf F^2$.
Its positive integral solutions are the squared Markov triples $(m^2, n^2, k^2)$, where $(m,n,k)$ is a Markov triple. 
These integral points are generated from the triple $(1,1,1)$ by means of the transformations
\begin{align*}
(x,y,z) \rightarrow (x,9xy-2x-2y-z,y), \quad (x,y,z) \rightarrow (y,9yz-2y-2z-x,z).
\end{align*}
and permutations of coordinates.
Our second task in this paper is to determine the composition law on $\mathbf F^2$ as well as the quadric group laws on the associated hyperplane sections.

\section{Algebraic geometry of the Fricke surface $\mathbf F$.}
Let us denote the projectivization of $\F$ inside $\mathbf P(\mathbb C^3)$ by $\overline \F$.
It has the equation $x^2s+y^2s+z^2s=3xyz$. 
It is easy to check that $[x,y,z,s]=[0,0,0,1]$ is the only singular point of $\overline \F$.
Besides the six finite lines 
$$
(t, \pm it, 0), \quad (0,t,\pm it) \quad (\pm it, 0, t),
$$
$\overline \F$ contains the lines  $[0,y,z,0]$, $[x,0,z,0]$ and $[x,y,0,0]$ at infinity. 
Consider the rational parametrization
\begin{align*}
\phi: \mathbb P^2(\C) &\to \overline\F \\
[p,q,r] &\to
\left[p(p^2+q^2+r^2): q(p^2+q^2+r^2):r(p^2+q^2+r^2):3pqr \right] 
\end{align*}
obtained via the pencil of lines through $[0,0,0,1]$. It has the inverse
\begin{align*}
\psi:[x,y,z,s]\in \overline\F\setminus \{[0:0:0:1]\}\to [x:y:z]\in \mathbb P^2(\C).
\end{align*}
The map $\phi$ is well defined except the points $[1:\pm i:0]$, $[0:1:\pm i]$, $[1:0:\pm i]$, 
which all lie on the quadric in $\mathbb P^2(\C)$ defined by the equation $p^2+q^2+r^2=0$, and $\phi$ sends the remaining points of this quadric to the point $[0:0:0:1]$.  Moreover, it sends the lines $p=0$, $q=0$, $r=0$ 
respectively to the lines $[0,y,z,0]$, $[x,0,z,0]$ and $[x,y,0,0]$; the lines lying at infinity of  
$\overline\F$. 

Note that $\phi$ is defined everywhere on $\mathbb P^2(\Q)$ and its inverse is well defined on  
 $\overline\F(\Q)^*:=\overline\F(\Q)\setminus\{(0,0,0)\}$.
Hence, there is a bijection
$$
\phi: \mathbb P^2(\Q) \to \overline\F(\Q)^*,
$$
such that the lines $p=0$, $q=0$, $r=0$ are sent to the lines at infinity of $\overline\F(\Q)^*$.
Any element of $\mathbb P^2(\Q)$ can be represented by a unique triple 
$[p:q:r]\neq (0,0,0)$ of integers, such that $\gcd(p,q,r)=1$. 
The corresponding element of $\overline\F(\Q)^*$ is 
$$
\phi([p,q,r])=\left[p(p^2+q^2+r^2): q(p^2+q^2+r^2):r(p^2+q^2+r^2):3pqr \right] \in \F(\Q)^*
$$
and all rational points can be obtained in this manner. These solutions are finite provided $pqr\neq 0$, in which case we can express them as
$$
\phi([p:q:r])=\left(\frac{p^2+q^2+r^2}{3qr}, \frac{p^2+q^2+r^2}{3pr}, \frac{p^2+q^2+r^2}{3pq}
\right)\in \F(\Q)^*
$$
By setting $p/r=P$, $q/r=Q$ we obtain a map 
\begin{align*}
(P,Q)\in (\Q^*)^2 \to \left(\frac{P^2+Q^2+1}{3Q}, \frac{P^2+Q^2+1}{3P}, \frac{P^2+Q^2+1}{3PQ}\right) \in \F(\Q)^*
\end{align*}
with the inverse
\begin{align*}
(x,y,z)\in \F(\Q)^* \to (x/z,y/z) \in (\Q^*)^2.
\end{align*}

\section{The composition law on $\F$.}
Suppose $P=(m,n,k)$ and $Q=(a,b,c)$ are two points of $\F$. The line $L$  through them has the parametrization
$$
x=(a-m)t+m, \quad y=(b-n)t+n, \quad z=(c-k)t+k
$$
The intersection $L\cap \F$ has the equation
$$
[(a-m)t+m]^2+[(b-n)t+n]^2+[(c-k)t+k]^2=3[(a-m)t+m][(b-n)t+n][(c-k)t+k]
$$
Since $t=0$ and $t=1$ are two solutions of this equation, the third solution is
\begin{align*}
t
=\frac{3(ank+bmk+cmn)-2(am+bn+ck)-3mnk}{3(a-m)(b-n)(c-k)}
\end{align*}
and the corresponding point $(x,y,z)$ can be found as below:
\begin{proposition}
If $P=(a,b,c)\in \F(K)$ and $Q=(m,n,k)\in \F(K)$ are two points of $\F$ defined over some field $K$ of characteristic 0, then the composition $P\circ Q=(x,y,z)\in \F(K)$ is given by
\begin{align}\label{composex}
x=\frac{3(ank+bcm)-2(am+bn+ck)}{3(b-n)(c-k)}\nonumber\\
y=\frac{3(bmk+acn)-2(am+bn+ck)}{3(a-m)(c-k)}\\
z=\frac{3(cmn+abk)-2(am+bn+ck)}{3(a-m)(b-n)}. \nonumber
\end{align}
In case  $(a,b,c)=(0,0,0)$ or $(a,b,c)=(m,n,k)$, the composition $(a,b,c)\circ (m,n,k)$ is not defined. 
Otherwise, if $(a-m)(b-n)(c-k)=0$, it is necessary to projectivize to get the answer
\begin{align*}
(m,n,k)\circ (m,b,c)=[m:n:k:1]\circ [m:b:c:1]=[0:b-n:c-k:0].
\end{align*}
Finally, the product is not well-defined if both of the points lies on the same line at infinity of $\F$. 
\end{proposition}

Note that the $\circ$-composition of two integral Markov triples is non-integral.

\medskip\noindent
{\bf Remark.} Let $\sigma$ be a constant. The similarly defined composition on the surface
$\F_\sigma: x^2+y^2+z^2=3xyz+\sigma$ has exactly the same expression
as (\ref{composex}). In case $\sigma=0$, one can bring 
(\ref{compose}) to the form: $P\circ Q=$
\begin{align*}
\left(\frac{(an-bm)^2+(ak-mc)^2}{3am(b-n)(c-k)},
\frac{(bk-cn)^2+(an-bm)^2}{3bn(a-m)(c-k)},
\frac{(ak-mc)^2+(bk-cn)^2}{3ck(a-m)(b-n)}\right)
\end{align*}
In order to carry out the analogy with the group law for elliptic curves, one may define the commutative law with identity
$$
(a,b,c)\star (m,n,k):=(1,1,1) \circ ((a,b,c)\circ (m,n,k)),
$$
whenever the right hand side is defined. 
However, $\star$ is not associative as one may easily check.

\subsection{Transfer of the structures to  $\mathbb P^2(\mathbb{Q})$}
We may transfer the Viete transformations and the composition $\circ$ to $\mathbb P^2$ by means of the parametrization $\phi$. The Viete transformations are given by 
\begin{align}\label{vietet}
&L: [p:q:r]\to [pr:p^2+q^2: qr]\\ 
&R: [p:q:r]\to [qp: q^2+r^2: pr]. \nonumber
\end{align}
Note that these are birational transformations which are not everywhere defined. 
On the other hand, viewed as transformations of 
$\mathbb P^2(\Q)$, they are well defined everywhere, except at the points $[0,0,1]$, $[0,1,0]$ and $[1,0,0]$.
We may apply permutations of coordinates to (\ref{vietet}) to get the involutions
\begin{align*}
&[p:q:r]\to [pr:qr: p^2+q^2]\\ {}
&[p:q:r]\to [q^2+r^2:pq: pr]\\ {}
&[p:q:r]\to [pq: p^2+r^2: rq],
\end{align*}
which defines a birational action of $\Z/2\Z*\Z/2\Z*\Z/2\Z$ on $\mathbb P^2$. We may also transfer the  composition $\circ$ to $\mathbb P^2(\Q)$ by means of $\phi$ as:
\begin{align*}
[a:b:c], [m:n:k] \in \mathbb P^2(\Q) \implies [a:b:c]\circ [m:n:k]=\\ {}
\Bigl[\left(  
( {a}^{2}+{b}^{2}+{c}^{2} ) kn- 
( {m}^{2}+{n}^{2}+{k}^{2}) bc 
\right)  
\left(  
(bm-an ) ^{2}+(cm-ak ) ^{2} 
\right),\\ 
 \left(  
( {a}^{2}+{b}^{2}+{c}^{2})km- 
( {m}^{2}+{n}^{2}+{k}^{2})a c 
 \right) 
 \left(  
( an-bm) ^{2}+ (cn-bk ) ^{2} 
 \right),\\ 
 \left( 
( {a}^{2}+{b}^{2}+{c}^{2} ) mn- 
( {m}^{2}+{n}^{2}+{k}^{2})ba 
 \right)  
 \left(  
( ak-cm ) ^{2}+ ( bk-cn ) ^{2} 
 \right).
 {}\Bigr]
\end{align*}

\subsection{Hyperplane sections of $\F$}
Let $H$ be a hyperplane in $\mathbb P^2(\C)$. 
Then $H\cap \overline\F$ is  a
a smooth cubic (if $H$ generic),
a nodal cubic (if $H$ passes through the origin or $H$ is a tangent to $\overline\F$) or
a union of three lines (if $H$ is the hyperplane at infinity) or
a quadric, together with a line at infinity (if $H$ is of the form $x=a$ or $y=b$ or $z=c$). In view of its connection to the Frobenius conjecture, we will be interested in the last case.

Suppose $(m_0,n_0,k_0)$ is a Markov triple. Then the plane $H:=\{(x,y,z): y=n_0\}$ contains infinitely many Markov triples. The intersection $H\cap \F$ is a quadric $Q_{n_0}$ with a special point ${(m_0,n_0,k_0)}$, given by the equation
\begin{align}\label{section}
Q_{n_0}: \quad x^2+n_0^2+z^2=3xn_0z
\end{align}
We may express $z$ explicitly in terms of $x$ as
\begin{align*}
z=\frac{3xn_0\pm\sqrt{9x^2n_0^2-4(n_0^2+x^2)}}{2}=
\frac{3xn_0\pm\sqrt{x^2n_0^2-4(x-n_0)^2)}}{2}
\end{align*}
Denote by $Q_{n_0}(\Z)$ the set of integral points and by $Q_{n_0}(\N)$ the set of positive integral points. We call a point of $Q_{n_0}(\N)$ with 
$n_0=\max\{m_0, n_0, k_0\}$ a {\it fundamental point}.  Note that the famous Frobenius conjecture is the statement that
every $Q_{n_0}(\N)$ has a unique fundamental point.

\paragraph{The infinity of $Q_{n_0}$.}
If we consider (\ref{section})
as a quadric in the $x-z$ plane, and projectivize it, we get the equation
$
x^2+n_0^2s^2+z^2=3xzn_0.
$
This quadric intersects the line $s=0$ at the points
$x^2+z^2=3xzn_0$. Deprojectivize by setting $t=x/z \implies$
$$
t^2-3tn_0+1=0 \implies t=\frac{3n_0\pm \sqrt{9n_0^2-4}}{2}.
$$
Then it is easy to check that
$$
\frac{3n_0+\sqrt{9n_0^2-4}}{2}=\lceil3n_0:3n_0:3n_0:3n_0:3n_0:\dots\rceil
$$
$$
\frac{3n_0-\sqrt{9n_0^2-4}}{2}=\lceil0:3n_0:3n_0:3n_0:3n_0:3n_0:\dots\rceil,
$$
where $\lceil n_0:n_1:\dots\dots\rceil=n_0-1/n_1-1/\dots$.
\subsubsection{The group law on the quadric sections $H\cap \F$.}
Let $O=(m_0,n_0,k_0)$ be a fundamental point.
Suppose $P_1:=(x_1, z_1)$ and $P_2:=(x_2, z_2)$ are on $Q_{n_0}$. 
If $P_1\neq P_2$ then the line through $O$ and parallel to $P_1P_2$ has the equation
$$
z=\frac{z_2-z_1}{x_2-x_1}(x-m_0)+k_0=\mu(x-m_0)+k_0
$$
and the intersection has the equation
$$
x^2+\left(\mu(x-m_0)+k_0\right)^2+n_0^2=3xn_0\left(\mu(x-m_0)+k_0\right) \iff
$$
Since $x=m_0$ is one solution of this equation the other solution is
\begin{align}\label{othersltn}
x=\frac{{n_0}^{2}+
 \left( -\mu m_0+k_0 \right) ^{2}}{1+\mu^2-3n_0 \mu}, \quad 
z=k_0-\mu\frac{2m_0+2\mu k_0-3n_0k_0-3m_0n_0\mu}{1+\mu^2-3n_0 \mu}
\end{align}
This gives the following result:
\begin{proposition}\label{qlaw}
Let $P_1:=(x_1, z_1)$ and $P_2:=(x_2, z_2)$ are on $Q_{n_0}$. Set 
$\mu:=(z_2-z_1)/(x_2-x_1)$.
If $P_1\neq P_2$, then $P_1\oplus P_2=(x,z)$, where 
$$
x=\frac{\mu^2m_0-m_0-2\mu k_0+3n_0k_0}{1+\mu^2-3n_0 \mu}, \quad z=\frac{k_0-2m_0\mu-\mu^2 k_0+3m_0n_0\mu^2}{1+\mu^2-3n_0 \mu}.
$$
If $P=(x_1, z_1)$, then $P\oplus P=(x,z)$, where
$$x=\frac{x_1^2n_0^2+(x_1k_0-z_1m_0)^2}{n_0^2m_0}, \quad
z=\frac{z_1^2n_0^2+(x_1k_0-z_1m_0)^2}{n_0^2k_0}. 
$$ 
\end{proposition}
It follows from this proposition that $(x_1,z_1)\oplus(z_1,x_1)=(k_0,m_0)$. 
Moreover, $(-m_0, -k_0)$ is an element of order 2 under $\oplus$.
\begin{proof}
The case $P_1\neq P_2$ is obtained from (\ref{othersltn}) by routine modifications.
It remains to establish the case $P_1=P_2$.
In this case, the line tangent to the quadric is given by the equation
$$
(2x_1-3n_0z_1)(x-x_1)+(2z_1-3n_0x_1)(z-z_1)=0,
$$
whose slope is 
$$
\sigma=-\frac{2x_1-3n_0z_1}{2z_1-3n_0x_1}
$$
The line through $(m_0,n_0,k_0)$ with this slope has the equation
$
z-k_0=\sigma(x-m_0),
$
and the second intersection of this line with $Q_{n_0}$ has the equation
$$
\implies x^2+n_0^2+\left(\sigma x-\sigma m_0+k_0\right)^2
=3xn_0\left(\sigma x- \sigma m_0+k_0\right).
$$
Since $m_0$ is a solution to this equation, the other solution is
$$
x=\frac{n_0^2+(-\sigma m_0+k_0)^2}{m_0(1+\sigma^2-3n_0\sigma)}.
$$
If we write $  z$ instead of $ x $ and $ m_0 $ instead of $ k_0, $ we get $ z. $
We obtain the desired result after some routine computations.
\end{proof}


\subsubsection{Inverses}
\begin{proposition}
The inverse of $P_1:=(x_1, z_1)$  is $(x,z)$ with
\begin{align*}
x&=\frac{(\mu^2x_1^2+z_1^2-2x_1\mu z_1+n_0^2)(2k_0-3n_0m_0)^2}{n_0^2(9n_0^2-4)} \\
z&=\mu{{ m_0}}^{2}-\mu{ x_1}+{ z_1}-\frac{4\mu}{9}+{\frac {\mu
 \left( \mu{ x_1}-{ z_1} \right) ^{2}{{ m_0}}^{2}}{{{ n_0}}^
{2}}}-{\frac {4\mu \left( 9 \left( \mu{ x_1}-{ z_1}
 \right) ^{2}+4 \right)}{9 \left(9{ n_0^2}-4 \right)   }}
\end{align*}
where $ \mu= -\frac{2m_0-3n_0k_0}{2k_0-3n_0m_0}.$
\end{proposition}
\begin{proof}
To find the inverse of $(x_1, z_1)$, we start with the tangent line at $O=(m_0,n_0,k_0)$:
$$
(2m_0-3n_0k_0)(x-m_0)+(2k_0-3n_0m_0)(z-k_0)=0
$$
Its slope is
$$
\mu=-\frac{2m_0-3n_0k_0}{2k_0-3n_0m_0}
$$
Hence the line through $(x_1,z_1)$ and parallel to the tangent line at $O$ has the equation
$
z=\mu(x-x_1)+z_1.
$
We are looking for the second intersection point of this line with the quadric
$$
x^2+n_0^2+\left(\mu(x-x_1)+z_1\right)^2=
3n_0x\left(\mu(x-x_1)+z_1\right).
$$
We know that one solution of this is $x_1$. Hence the other solution is
$$
\frac{\mu^2x_1^2+z_1^2-2x_1\mu z_1+n_0^2}{x_1(1+\mu^2-3n_0\mu)}.
$$
One has
$$
1+\mu^2-3n_0\mu=\frac{n_0^2(9n_0^2-4)}{(2k_0-3n_0m_0)^2},
$$
and we obtain desired result after some routine computations.
\end{proof}

\subsubsection{Connections between the Viete transformations and the group law}
Note that, if $(m,k)\in Q_{n_0}(\Z)$, then so are
\begin{align*}
\mathsf A(m,k):=(m, 3mn_0-k),\quad \mathsf T\mathsf A(m,k):=(3mn_0-k,m),   \\ \mathsf C(m,k):=(3n_0k-m,k), \quad \mathsf T\mathsf C(m,k):=(k, 3n_0k-m),\\
\mathsf B(m,k)=(-m,-k)
\end{align*}
where $\mathsf T(m,k):=(k,m)$ is the transposition. Since $\mathsf A^2=\mathsf T^2=Id$ and $\mathsf C=\mathsf T\mathsf A\mathsf T$, the transformations $\mathsf A$ and $\mathsf T$ generate the infinite dihedral group.
Hence, the Markov triples lying on $Q_{n_0}$ are produced by Viète transformations, starting from a fundamental point 
$(m_0, n_0, k_0)$.

Suppose $P:=(x_1, z_1)\in Q_{n_0}$. Then 
$\mathsf C P=(3n_0z_1-x_1, z_1)$. In the computation of $P\oplus \mathsf CP$, one has  $\mu=0$ in Proposition~\ref{qlaw}, so that 
$$
P\oplus \mathsf CP =(3n_0k_0-m_0,k_0)=\mathsf CO
$$
Similarly, $\mathsf AP=(x_1,3n_0x_1-z_1)$ and in the computation of $P\oplus \mathsf AP$, one has  $\mu=\infty$ in Proposition~\ref{qlaw}, so that 
$$
P\oplus \mathsf AP =(m_0,3n_0m_0-k_0)=\mathsf AO
$$


The transformation
$$
\mathsf T \mathsf A:(m,n_0,k)\to(3mn_0-k,n_0, m)
$$
can be seen as a linear transformation $(m,k)\to(3mn_0-k,m)$, i.e.
$$
\left(\begin{matrix}
3n_0 & -1\\ 1&0 
\end{matrix}\right)\left(
\begin{matrix}
m\\k
\end{matrix}\right)
=
\left(\begin{matrix}
3mn_0-k,&m
\end{matrix}\right)
$$
One has 
$$
\left(\begin{matrix}
3n_0 & -1\\ 1&0 
\end{matrix}\right)^r
=
\left(\begin{matrix}
b_r & -b_{r-1}\\ b_{r-1}&-b_{r-2} 
\end{matrix}\right)
$$
where $b_r$ is the Chebyshev-like polynomial defined by the recursion
$$
b_{r+2}(n_0)=3n_0b_{r+1}(n_0)-b_r, \quad b_0(n_0)=1, \quad b_1(n_0)=3n_0
$$
then
$$
 b_2(n_0)=9n_0^2-1, \quad b_3(n_0)=27n_0^3-6n_0, \quad b_5(n_0)=81n_0^4-27n_0^2+1,...
$$
so that 
$$
\frac{b_{r+2}(n_0)}{b_{r+1}(n_0)}=3n_0-\frac{1}{\frac{b_{r+1}(n_0)}{b_{r}(n_0)}}. 
$$
Hence
$$
\left(\begin{matrix}
b_r & -b_{r-1}\\ b_{r-1}&-b_{r-2} 
\end{matrix}\right)
\left(
\begin{matrix}
m\\k
\end{matrix}\right)
=(mb_r(n_0)-kb_{r-1}(n_0), mb_{r-1}(n_0)-kb_{r-2}(n_0)).
$$
So that 
\begin{eqnarray*}
\mathsf T \mathsf A^r(m,n_0,k)=(mb_r(n_0)-kb_{r-1}(n_0),n_0, mb_{r-1}(n_0)-kP_{r-2}(n_0))\\
\mathsf T \mathsf C^r(m,n_0,k)=(kb_r(n_0)-mb_{r-1}(n_0),n_0, kb_{r-1}(n_0)-mb_{r-2}(n_0)).
\end{eqnarray*}
This yields the following lemma.
\begin{lemma} \label{lemma}
Let $ {(\mathsf T \mathsf A)}^r =\mathsf T \mathsf A{(\mathsf T \mathsf A)}^{r-1}$ and  $ {(\mathsf T \mathsf C)}^r =\mathsf T \mathsf C{(\mathsf T \mathsf C)^{r-1}} $ denote respectively $ r- $times composition of $\mathsf T \mathsf A. $ and $\mathsf T \mathsf C.$ Then it is given by 
\begin{equation*}
{\mathsf T \mathsf A}^r:(m,n_0,k) \to (mb_r(n_0)-kb_{r-1}(n_0),n_0,mb_{r-1}(n_0)-kb_{r-2}(n_0)).
\end{equation*}
\begin{equation*}
{\mathsf T \mathsf C}^r:(m,n_0,k) \to (kb_{r}(n_0)-mb_{r-1}(n_0),n_0,kb_{r-1}(n_0)-mb_{r-2}(n_0)).
\end{equation*}
\end{lemma}
We have shown with Maple implementations and  Lemma \ref{lemma} that elements in $Q_{n_0}$  formed with operators $\mathsf T \mathsf A$  and $\mathsf T \mathsf C $  always contain an integral subgroup.
It can be proved by induction. 
Firstly we will show that the line parallel to  $\mathsf T \mathsf A$ and the lines passing through points $O(m_0,n_0,k_0)$ and ${\mathsf T \mathsf A}^2$ are parallel. Clearly it is sufficient to show that their slopes are equal.
For $ r=2 $  slope of the tangent line at the point $\mathsf T \mathsf A  $ is equal slope of the line passing through $O(m_0,n_0,k_0)$ and ${\mathsf T \mathsf A}^2. $
Equation of the tangent line at  $\mathsf T \mathsf A  $ 
\begin{equation} \label{TA}
(2(3mn_0-k)-3n_0m)(x-3mn_0+k)+(2m-3(3mn_0-k)n_0)(z-m)=0
\end{equation}
hence slope of the tangent line \eqref{TA} 
$$ \mu_1=\dfrac{2k-3mn_0}{2m-9n_0^2m+3n_0k}. $$
On the other hand equation of the line passing through $O(m_0,n_0,k_0)$ and ${\mathsf T \mathsf A}^2 $
\begin{equation} \label{0andT2}
(2k-3n_0m)(x-m)=(z-k)(2m-9n_0^2m+3n_0k+m)
\end{equation}
then slope of the equation \eqref{0andT2} 
$$ \mu_2=\dfrac{2k-3mn_0}{2m-9n_0^2m+3n_0k}. $$ Clearly $ \mu_1=\mu_2. $
For $ r=n-1$ we will assume that slope of line passing through $O(m_0,n_0,k_0)$ and ${\mathsf T \mathsf A}^{n-1} $ is equal slope of the tangent line passing through  ${\mathsf T \mathsf A}^{n-2} $ and ${\mathsf T \mathsf A}^{n-1}. $  Using the induction hyphothesis and simple calculations we proved
slope of line passing through $O(m_0,n_0,k_0)$ and ${\mathsf T \mathsf A}^{n} $ is equal slope of the tangent line passing through  ${\mathsf T \mathsf A}^{n-1} $ and ${\mathsf T \mathsf A}^{n}. $ The similar process can be repeated for ${\mathsf T \mathsf C}. $ As a result elements in $Q_{n_0}$  formed with operators $\mathsf T \mathsf A$  and $\mathsf T \mathsf C $  always contain an integral subgroup.

\section{The double Fricke surface $\F^2$}
Recall that $\F^2$ is the surface given by the equation $(x+y+z)^2=9xyz$.
Let us denote the projectivization of $\F^2$ inside $\mathbf P(\mathbb C^3)$ by $\overline \F^2$.
It has the equation $s(x+y+z)^2=9xyz$. 
It is easy to check that $[x,y,z,s]=[0,0,0,1]$ is the only singular point of 
$\overline \F^2$.
Besides the three finite lines 
$$
(t, -t, 0), \quad (0,t,-t) \quad ( t, 0, -t),
$$
$\overline \F^2$ contains the lines  $[0,y,z,0]$, $[x,0,z,0]$ and $[x,y,0,0]$ at infinity. 
The rational points of this surface has the parametrization
\begin{align*}
(P,Q)\in (\Q^*)^2 \to \left(\frac{(P^2+Q^2+1)^2}{9Q^2}, \frac{(P^2+Q^2+1)^2}{9P^2}, \frac{(P^2+Q^2+1)^2}{9P^2Q^2}\right) \in \F^2(\Q^*)
\end{align*}
as one may check.
\begin{figure}[ht]
\begin{center}
\includegraphics[width=6.5cm]{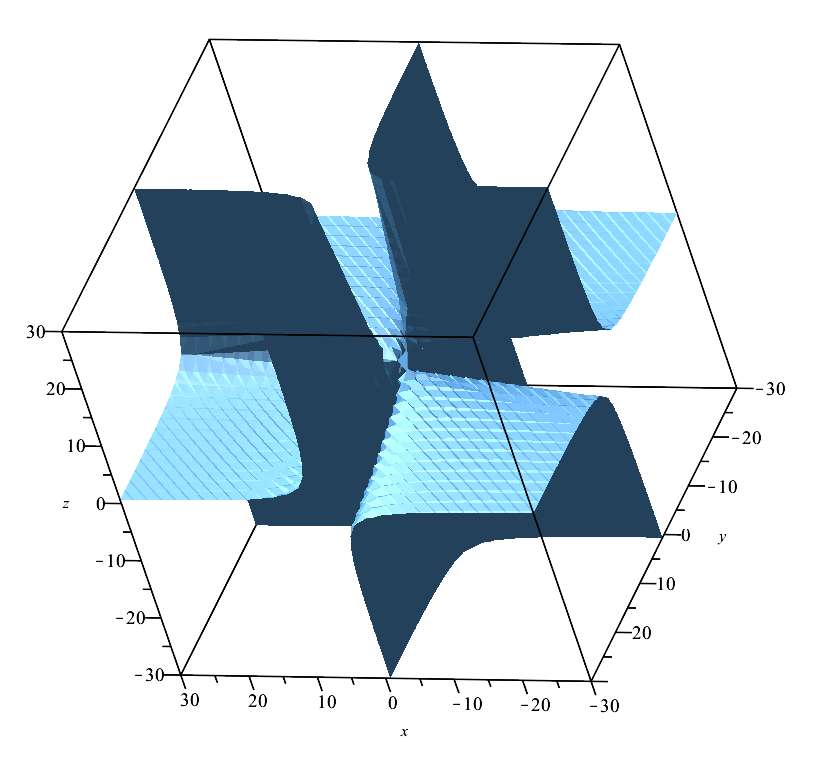}\quad
\includegraphics[width=6.5cm]{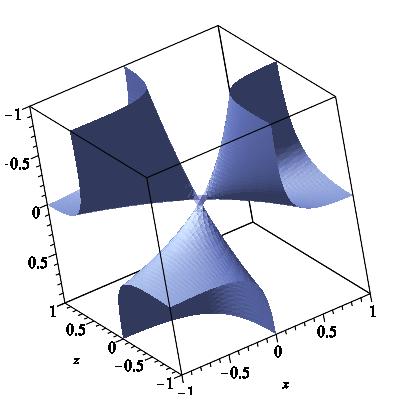}
\caption{Fricke surface (left) and the double Fricke surface}
\end{center}
\end{figure}
\begin{theorem}
Positive integral points of $\F^2$ are precisely those of the form $(m^2,n^2,k^2)$, where $(m,n,k)$ is a Markov triple. All positive integral solutions are obtained from the solution triple $(1,1,1)$ by use of the transformations
\begin{align}\label{nielsen}
(x,y,z) \rightarrow (x,9xy-2x-2y-z,y), \quad (x,y,z) \rightarrow (y,9yz-2y-2z-x,z).
\end{align}
and permutations of coordinates.
\end{theorem}
\begin{proof}
Suppose that $(m,n,k)$ is a Markov triple. Then 
$$
m^2+n^2+k^2=3mnk\implies (m^2+n^2+k^2)^2=9m^2n^2k^2,
$$
i.e. $(m^2, n^2, k^2)$ is on $\F^2(\Z_{>0})$. 
Now suppose $(x,y,z)\in \F^2(\Z_{>0})$. Then 
$(\sqrt x, \sqrt y, \sqrt z)\in \F$, since
$$
(x+y+z)^2=9xyz\implies x+y+z=3\sqrt{xyz}\implies
(\sqrt x)^2+ (\sqrt y)^2+ (\sqrt z)^2=3\sqrt{xyz}
$$
 Let $ (m,n,k) $ is solution triple. Then $ m $ is a root of the polynomial 
$$f(x)=(x+n+k)^2-9xnk=x^2+(2n+2k-9nk)x+n^2+k^2+2nk.  $$ The other root $m_0$ thus satisfies
\begin{equation} \label{eq.}
m_0=9nk-2n-2k-m=\frac{(n+k)^2}{m}.
\end{equation}
The first equation in \eqref{eq.} says that $ m_0 $ is integer and second that $ m_0 $ is positive. If $ (m,n,k)=(1,1,1) $ then $ (m_0,1,1) $ is another solution triple. Repeating with second and third element of triple, we see that for every triple $(m,n,k)  $ we get three others.
\end{proof}
{\bf Claim.}  Solution triples of $(9yz-2y-2z-x,y,z)$ \eqref{nielsen} obtained from the fundamental solution $(1,1,1)$  are the triples $(x^2,y^2,z^2)$, where $(x,y,z)$ is a Markov triple.
Solution triples of Markov equation is $(3yz-x,y,z)$.
\begin{align*}
((3yz-x)^2,y^2,z^2)=(9y^2z^2-6xyz+x^2,y^2,z^2)
\end{align*}
We know that $ (x,y,z) $ is Markov triple then we can get
$$ (9y^2z^2-2y^2-2z^2-x^2,y^2,z^2). $$
Note that $\F^2$ has other integral points.
For example, we have the following solutions produced from the ``{fundamental solution}" $(0,1,-1)$:
\begin{align*}
(0, 1, -1),(-1, -9, 1), (1, -64, -9),(-9, 100, -1) ,
 (100, -8281, -9), \dots 
\end{align*}
\begin{proposition}
Let $(a,b,c)\in \F^2(\Z)$ such that $a$, $b$, $c$ are not simultaneously positive. Then $(a,b,c)$ is produced from a ``{fundamental solution}" $(-n,0,n)$ with $n\in \Z$ by use of the transformations (\ref{nielsen}) and its permutations.  
\end{proposition}
\begin{proof}
 Let $(a,b,c)\in \F^2(\Z)$. Then $ a $ is a root of the polynomial 
$$f(x)=(x+b+c)^2-9xbc=x^2+(2b+2c-9bc)x+b^2+c^2+2bc.  $$ The other root $a_0$ thus satisfies
\begin{equation} \label{eq.}
a_0=9bc-2b-2c-a=\frac{(b+c)^2}{a}.
\end{equation}
Hence $(a_0,b,c)  $ is another Markov triple. Repeating with $b  $ and $ c $, we see that for every nonsingular
triple $ (a,b,c) $ we get three others.
 \begin{align}\label{transf}
 (a_0=9bc-2b-2c-a,b,c)  \\  \nonumber
(a,b_0=9ac-2a-2c-b,c)  \\ \nonumber
(a,b,c_0=9ab-2a-2b-c)  \nonumber
 \end{align}
Suppose $ b=0 $ is a root of the polynomial then we get from \eqref{transf} and double Markov equation $ a=-c $. 
On the other hand suppose without loss of the generality, $a=-c $. Since $ (a,b,c)\in \F^2(\Z) $ then $ b^2=-9a^2b $. This gives $ b=0. $ 
\end{proof}

Computer experiments indicate that the Frobenius conjecture also holds for the solution tree generated from the fundamental solution $(0,n,-n)$, in the sense that if $(a,b,c)$ belongs to this tree, then the maximum of the numbers $|a|, |b|, |c|$ determines the remaining two.

\medskip\noindent
{\bf Remark.} Equation~\ref{gen}  appears in Perrine~\cite{perrine} (Pg 229, Théorème 2.1) as a member of a special family (E), though it has not been given a special treatment. 
More generally, positive integral solution $n$-tuples of the equation 
$(x_1+\dots+x_n)^2=k^2x_1\dots x_n$ 
are the squares of the solution $n$-tuples of the {\it Hurwitz equation}
$
x_1^2+\dots +x_n^2=kx_1\dots x_n.
$
\subsection{The composition law on $\F^2$.}
Suppose $(m,n,k)$ and $(a,b,c)$ are two rational squared Markov triples. 
Let $L$ be the line through them. It has the parametrization
$$
x=(a-m)t+m, \quad y=(b-n)t+n, \quad z=(c-k)t+k
$$
Its intersection with the Fricke surface has the equation
$$
[(a-m)t+m+(b-n)t+n+(c-k)t+k]^2=9[(a-m)t+m][(b-n)t+n][(c-k)t+k]
$$
 Since $ t=0 $ and $ t=1 $ are two solutions of this equation, the third solution is 
 \begin{align*}
t=\frac{9(a-m)nk+9(b-n)mk+9(c-k)mn-2(a+b+c-m-n-k)(m+n+k)}{9(a-m)(b-n)(c-k)}
\end{align*}
Hence the third intersection point has the coordinates
$$
x=(a-m)t_0+m, \quad y=(b-n)t_0+n, \quad z=(c-k)t_0+k.
$$
Routine algebraic manipulations yield the following result:
\begin{proposition}
If $P=(a,b,c)\in \F^2(K)$ and $Q=(m,n,k)\in \F^2(K)$ are two points of $\F^2$ defined over some field $K$ of characteristic 0, then the composition $P\circ Q=(x,y,z)\in \F^2(K)$ is given by
\begin{align}\label{compose}
x=\frac{9(ank+bcm)-2(a+b+c)(m+n+k)}{9(b-n)(c-k)}\nonumber\\
y=\frac{9(bmk+acn)-2(a+b+c)(m+n+k)}{9(a-m)(c-k)}\\
z=\frac{9(cmn+abk)-2(a+b+c)(m+n+k)}{9(a-m)(b-n)} \nonumber
\end{align}
In case  $(a,b,c)=(0,0,0)$ or $(a,b,c)=(m,n,k)$, the composition $(a,b,c)\circ (m,n,k)$ is not defined. 
Otherwise, if $(a-m)(b-n)(c-k)=0$, it is necessary to projectivize to get the answer
\begin{align*}
(m,n,k)\circ (m,b,c)=[m:n:k:1]\circ [m:b:c:1]=[0:b-n:c-k:0].
\end{align*}
Finally, the product is not well-defined if both of the points lies on the same line at infinity of $\F^2$. 
\end{proposition}

Consequently, we define the operation
$$
(a,b,c)\overline\circ (m,n,k):=(x_0, y_0, z_0)
$$
This is commutative though not associative, and it satisfies 
\begin{align*}
\alpha \overline\circ \beta=\gamma \iff \alpha \overline\circ \gamma=\beta \iff \beta \overline\circ \gamma=\iff \alpha
\end{align*}
Note that the $\overline\circ$-product of two integral squared Markov triples is usually a non-integral rational squared Markov triple. Also note that this is not the same product on $\F(\Q)$, i.e. 
if $(a,b,c), (m,n,k)\in \F(\Q)$ and $(a,b,c)\circ (m,n,k)=(p,q,r)$ then
$$
(a^2,b^2,c^2)\overline\circ (m^2,n^2,k^2) \neq (p^2,q^2, r^2)
$$
 Let $ (a,b,c)=(2,1,1) $ and $ (m,n,k)=(1,2,5) $ then $$(p,q,r)=(2,1,1)\circ (1,2,5)=\left(\frac{15}{4},-\frac{3}{4},-6\right). $$
If $ (a^2,b^2,c^2)=(4,1,1) $ and $ (m^2,n^2,k^2)=(1,4,25). $ then
$$ (4,1,1) \overline\circ (1,4,25)=\left(\frac{361}{72},-\frac{7}{24},-\dfrac{28}{3}\right)$$ so we get $$
(a^2,b^2,c^2)\overline\circ (m^2,n^2,k^2) \neq (p^2,q^2, r^2).
$$
{
\subsection{Transfer of the structures to  $\mathbb P^2(\mathbb{Q})$}
We may transfer the Viete transformations and the composition $\circ$ to $\mathbb P^2(\mathbb{Q})$ by means of the parametrization $\phi$. The Viete transformations are given by 
\begin{align}\label{vietet}
&L: [p:q:r]\to [pr:(p+q)^2: qr]\\ 
&R: [p:q:r]\to [qp: (q+r)^2: pr]. \nonumber
\end{align}
Note that these are birational transformations which are not everywhere defined. 
On the other hand, viewed as transformations of 
$\mathbb P^2(\Q)$, they are well defined everywhere, except at the points $[0,0,1]$, $[0,1,0]$ and $[1,0,0]$.
We may apply permutations of coordinates to (\ref{vietet}) to get the involutions
\begin{align*}
&[p:q:r]\to [pr:qr: (p+q)^2]\\ {}
&[p:q:r]\to [(q+r)^2:pq: pr]\\ {}
&[p:q:r]\to [pq: (p+r)^2: rq],
\end{align*}
which defines a birational action of $\Z/2\Z*\Z/2\Z*\Z/2\Z$ on $\mathbb P^2$. We may also transfer the  composition $\circ$ to $\mathbb P^2(\Q)$ by means of $\phi$ as: 
\begin{align*}
[a:b:c], [m:n:k] \in \mathbb P^2(\Q) \implies [a:b:c]\circ [m:n:k]=\\ {}
\Bigl[\left(  
( {a}+{b}+{c} )^{2} kn- 
( {m}+{n}+{k})^{2} bc 
\right)  
\left(  
(ank+bcm ) -2(a+b+c )(m+n+k) 
\right),\\ 
 \left(  
( {a}+{b}+{c} )^{2}km- 
( {m}+{n}+{k})^{2}a c 
 \right) 
 \left(  
(bmk+acn ) -2(a+b+c )(m+n+k) 
 \right),\\ 
 \left( 
( {a}+{b}+{c} )^{2} mn- 
( {m}+{n}+{k})^{2}ba 
 \right)  
 \left(  
(cmn+abk ) -2(a+b+c )(m+n+k)
 \right).
 {}\Bigr]
\end{align*}
}
\subsection{Quadric sections of $\F^{2}$}
Suppose $(m_0,n_0,k_0)$ is a squared Markov triple. Then the plane $H:=\{(x,y,z): y=n_0\}$ contains infinitely many squared Markov triples. The intersection $H\cap \F^2$ is a quadric $Q_{n_0}\subset H$ with a special point ${(m_0,n_0,k_0)}$, given by the equation
\begin{align}
Q_{n_0}: \quad (x+n_0+z)^2=9xn_0z
\end{align}
Consider this as a quadric in the $x-z$ plane, and projectivize by $s$ to get the equation
$$
(x+n_0s+z)^2=9xzn_0
$$
This quadric intersects the line $s=0$ at the points
$(x+z)^2=9xzn_0$. Deprojectivize by setting $t=x/z$, we find the points at infinity of $Q_{n_0}$ to be
$$
t=\frac{3n_0-2\pm \sqrt{(9n_0-2)^2-4}}{2}
$$

\subsubsection{The group law on the quadric sections $H\cap \F^2$.}


\begin{proposition}
Let $P_1:=(x_1, z_1)$ and $P_2:=(x_2, z_2)$ are on $Q_{n_0}$. Set 
$\sigma=-(2x_1+2n_0+2z_1-9n_0z_1)/(2z_1+2n_0+2x_1-9n_0x_1).$
If $P_1\neq P_2$, then $P_1\oplus P_2=(x,z)$, where 
\begin{align*}
x=&\frac{9n_0k_0(x_2-x_1)^2-(2n_0+2k_0)(x_2-x_1)(z_2-z_1)}{(z_2-z_1+x_2-x_1)^2-9n_0(x_2-x_1)(z_2-z_1)} \\
&-\frac{(m_0+2n_0+2k_0)(x_2-x_1)^2+m_0(z_2-z_1)^2}{(z_2-z_1+x_2-x_1)^2-9n_0(x_2-x_1)(z_2-z_1)},
\end{align*}
\begin{align*}
z=\frac{(z_2-z_1)^2(9n_0m_0-2m_0-2n_0-k_0)-2(z_2-z_1)(x_2-x_1)(m_0+n_0)+k_0}{-2n_0(n_0+x_2+z_2+x_1+z_1)+(9n_0-2)(z_2x_1+z_1x_2)-2(z_1z_2+x_1x_2)}.
\end{align*}
If $P=(x_1, z_1)$, then $P\oplus P=(x,z)$, where
$$
x=\frac{\left[9n_0(n_0+k_0)x_1+m_0z_1]^2-4(n_0+k_0+m_0)[x_1(x_1+n_0)(n_0+k_0)+z_1(z_1+n_0)m_0\right]}{9n_0^3m_0},
$$
\begin{align*}
z=\frac{\left[9n_0(n_0+m_0)z_1+k_0x_1]^2-4(n_0+m_0+k_0)[z_1(z_1+n_0)(n_0+m_0)+x_1(x_1+n_0)k_0\right]}{9n_0^3k_0}.
\end{align*}
\end{proposition}
\begin{proof}
We take $O=(m_0,n_0,k_0)$ as the neutral element. Suppose $P_1:=(x_1, z_1)$ and $P_2:=(x_2, z_2)$ are on $Q_{n_0}$. 
If $P_1\neq P_2$ then the line through $O$ and parallel to $P_1P_2$ has the equation
$$
z=\frac{z_2-z_1}{x_2-x_1}(x-m_0)+k_0
$$
and the intersection has the equation
$$
\left(x+\left(\frac{z_2-z_1}{x_2-x_1}(x-m_0)+k_0\right)+n_0\right)^2=9xn_0\left(\frac{z_2-z_1}{x_2-x_1}(x-m_0)+k_0\right)
$$
Set $u:=x-m_0$ and $\mu={z_2-z_1}/{x_2-x_1}$.  Then the equation becomes 
$$
\left(u+m_0+\left(\mu u+k_0\right)+n_0\right)^2=9(u+m_0)n_0\left(\mu u+k_0\right)
$$
and we know that $x=m_0 \iff u=0$ is one solution of this equation. 
\begin{align*}
\left(u(\mu+1)+m_0+n_0+k_0\right)^2=9(u+m_0)n_0\left(\mu u+k_0\right)\\
u=\frac{-2(\mu+1)(m_0+n_0+k_0)
+9n_0k_0+9n_0m_0 \mu}{(\mu+1)^2-9n_0\mu} 
\end{align*}
Thus $x=u+m_0$, $z=u\mu+k_0$ $\implies$
\begin{align*}
x=\frac{9n_0k_0-2\mu n_0-2\mu k_0-m_0-2n_0-2k_0+m_0{\mu}^2}{(\mu+1)^2-9n_0\mu} 
\end{align*}
\begin{align*}
z=\frac{9n_0m_0{\mu}^2-2{\mu}^2m_0-2{\mu}^2n_0-{\mu}^2k_0-2\mu m_0-2\mu n_0+k_0}{(\mu+1)^2-9n_0\mu} 
\end{align*}
If $ z_1=z_2=k $ then $ \mu=0 $ and $ x=9n_0k_0-2n_0-2k_0-m_0 $ $ z=k_0 $. Hence
\begin{align*}
(x_1,k)\circ(x_2,k)=(9n_0k_0-2n_0-2k_0-m_0,k_0)
\end{align*}
Note that $ z_1=z_2=k $ happens when $ (x_2,z_2)=(9n_0z_1-2n_0-2z_1-x_1,z_1). $ Similarly, if $ x_1=x_2=k $ then $ \mu=\infty $ and $ x=m_0 $, $ z=9n_0m_0-2m_0-2n_0-k_0$. Hence 
\begin{align*}
(m,z_1)\circ(m,z_2)=(m_0,9n_0m_0-2m_0-2n_0-k_0).
\end{align*}
Note that $ x_1=x_2=m $ happens when $ (x_2,z_2)=(x_1,9n_0x_1-2n_0-2x_1-z_1). $

For general $ \mu $ we have 
\begin{align*}
x=&\frac{9n_0k_0(x_2-x_1)^2-(2n_0+2k_0)(x_2-x_1)(z_2-z_1)}{(z_2-z_1+x_2-x_1)^2-9n_0(x_2-x_1)(z_2-z_1)} \\
&-\frac{(m_0+2n_0+2k_0)(x_2-x_1)^2+m_0(z_2-z_1)^2}{(z_2-z_1+x_2-x_1)^2-9n_0(x_2-x_1)(z_2-z_1)}
\end{align*}
For the denominator one has
\begin{align*}
(z_2-z_1)^2+2(z_2-z_1)(x_2-x_1)+(x_2-x_1)^2-9n_0(x_2-x_1)(z_2-z_1) \\
=-2n_0(n_0+x_2+z_2+x_1+z_1)+(9n_0-2)(z_2x_1+z_1x_2)-2(z_1z_2+x_1x_2)
\end{align*}
and
\begin{align*}
z=\frac{(z_2-z_1)^2(9n_0m_0-2m_0-2n_0-k_0)-2(z_2-z_1)(x_2-x_1)(m_0+n_0)+k_0}{-2n_0(n_0+x_2+z_2+x_1+z_1)+(9n_0-2)(z_2x_1+z_1x_2)-2(z_1z_2+x_1x_2)}.
\end{align*}
It remains to establish the case $ P=Q. $ In this case, the line tangent to the quadric is given by the equation
$$
(2x_1+2n_0+2z_1-9n_0z_1)(x-x_1)+(2z_1+2n_0+2x_1-9n_0z_1)(z-z_1)=0,
$$
whose slope is 
$$
\sigma=-\frac{2x_1+2n_0+2z_1-9n_0z_1}{2z_1+2n_0+2x_1-9n_0x_1}.
$$
The line through $(m_0,n_0,k_0)$ with this slope has the equation
$$
z-k_0=-\frac{2x_1+2n_0+2z_1-9n_0z_1}{2z_1+2n_0+2x_1-9n_0x_1}(x-m_0),
$$
and the second intersection of this line with $Q_{n_0}$ has the equation
\begin{align*}
\left(x+n_0+\left(\sigma(x-m_0)+k_0\right)\right)^2 
=9xn_0\left(\sigma(x-m_0)+k_0\right).
\end{align*}
We know that $m_0$ is a solution to this equation.
Hence the second solution is
$$
x=\frac{\left[(n_0+k_0)(2z_1+2n_0+2x_1-9n_0x_1)+m_0(2x_1+2n_0+2z_1-9n_0z_1)\right]^2}{m_0[(9n_0z_1-9n_0x_1)^2+9n_0(2z_1+2n_0+2x_1-9n_0z_1)(2z_1+2n_0+2x_1-9n_0x_1)]}
$$
We obtain desired result after some routine computations and simplification.
\end{proof}

\subsubsection{Inverses}
\begin{proposition}
The inverse of $P_1:=(x_1, z_1)$  is $(x,z)$ with
\begin{align*}
x&=\frac{{\sigma}^{2}{{x_1}}^{2}-2\sigma{x_1}{z_1}+{{n_0}}^{2}+2{ n_0}{ z_1}-2{
n_0}\sigma{ x_1}+{{ z_1}}^{2}
}{81n_0^4m_0} \\
z&=-\sigma{x_1}+{ z_1}+{\frac {2\sigma\left( {
z_1}+\sigma{x_1} \right) }{{81{ n_0}}^{3}{m_0}}}+{\frac {\sigma}{{81{ n_0}}^{2}{m_0}}}-{
\frac {\sigma\left( {\sigma}{{x_1}}-{{z_1}} \right)^2 }{{81{n_0}}^{4}{ m_0}}} 
\end{align*}
where $ \sigma=-\frac{2x_1+2n_0+2z_1-9n_0z_1}{2z_1+2n_0+2x_1-9n_0x_1}.$
\end{proposition}
\begin{proof}
To find the inverse of $(x_1, z_1)$, we start with the tangent line at $O$:
$$
(2x_1+2n_0+2z_1-9n_0z_1)(x-m_0)+(2z_1+2n_0+2x_1-9n_0x_1)(z-k_0)=0
$$
Its slope is
$$
\sigma=-\frac{2x_1+2n_0+2z_1-9n_0z_1}{2z_1+2n_0+2x_1-9n_0x_1}=-\frac{m_0'-m_0}{k_0'-k_0}.
$$
Hence the line through $(x_1,z_1)$ and parallel to the tangent line at $O$ has the equation
$
z=\sigma(x-x_1)+z_1.
$
We are looking for the second intersection point of this with the quadric:
$$
\left(x+n_0+\sigma(x-x_1)+z_1\right)^2=9xn_0(\sigma(x-x_1)+z_1)
$$
We know that one solution of this is $x_1$. Hence the other solution is
$$
\frac{{\sigma}^{2}{{x_1}}^{2}-2\,\sigma{ x_1}\,{z_1}+{{n_0}}^{2}+2{n_0}{z_1}-2{n_0}\sigma{x_1}+{{z_1}}^{2}
}{ 2\sigma+1-9\,{n_0}\sigma+{\sigma}^{2} }.
$$
 Note that the denominator is independent of $(x_1, z_1)$.
$$
=(9n_0)^2\left(-x_1^2-z_1^2-2x_1z_1-2n_0z_1-2n_0x_1+9n_0x_1z_1\right)=81n_0^4m_0
$$
\end{proof}

\medskip\noindent
{\bf Acknowledgements}
This work is supported by the T\"{U}B\.{I}TAK grant 115F412 and the Galatasaray University research grant ****.

\end{document}